\documentclass[a4paper,11pt]{article}

\usepackage{authblk}

\usepackage{amssymb}
\let\Emptyset\emptyset

\usepackage{type1cm}        
\usepackage{makeidx}         
\usepackage{graphicx}        
\usepackage{multicol}        

\usepackage{newtxtext}       %
\usepackage{newtxmath}       

\makeindex             

\usepackage{mdframed}

\usepackage{tikz}
\usetikzlibrary{shapes,arrows}
\usetikzlibrary{patterns}

\usepackage{pgfplots}
\pgfplotsset{compat=1.13}
\usepgfplotslibrary{fillbetween}

\definecolor{HSUred}{RGB}{197,0,66}

\usepackage{hyperref}
\usepackage{xcolor}
\usepackage{subfig}
\usepackage{graphicx}

\renewcommand{\emptyset}{\Emptyset}

\newcommand{\doi}[1]{\href{https://doi.org/#1}{doi:#1}}

\newtheorem{defi}{Definition}[section]
\newtheorem{theorem}[defi]{Theorem}

\newtheorem{remark}[defi]{Remark}

\newenvironment{mproof}{\paragraph{Proof.}}{\hfill$\blacksquare$}

\numberwithin{equation}{section}
\numberwithin{table}{section}
\numberwithin{figure}{section}

\usepackage{geometry}
\begin{document}

\title{\Large \textbf{Numerical Study of Approximation Techniques for the Temporal Weights
to the DWR Method}}
\author[M.\ P.\ Bruchh\"auser, M.\ Bause]
{\large \textbf{Marius Paul Bruchh\"auser}\thanks{bruchhaeuser@hsu-hamburg.de ($^\ast$corresponding author)} 
$\,\boldsymbol{\cdot}$
\textbf{Markus Bause}\thanks{bause@hsu-hamburg.de}\\
{\small Helmut Schmidt University, University of the German Federal Armed Forces Hamburg,\\
Faculty of Mechanical and Civil Engineering, Chair of Numerical Mathematics\\ 
Holstenhofweg 85, 22043 Hamburg, Germany}
}
\date{}
\maketitle

\begin{abstract}
\noindent
This work presents a numerical investigation of different approximation 
techniques for the temporal weights used in the Dual Weighted Residual (DWR) 
method applied to a time-dependent convection-diffusion equation which is assumed 
to be convection-dominated.  
It is a continuation of a previous work by the authors where spatial weights
were compared for a steady-state case. 
A higher-order finite elements approach is compared to a more cost-efficient
higher-order reconstruction approach. 
Numerical examples point out the results regarding accuracy, efficiency
and stability reasons.
\end{abstract}


\section{Introduction}
\label{sec:1:introduction}
The DWR method has attracted researchers' interest in 
many fields of application problems since it was introduced by Becker and 
Rannacher \cite{BruchhaeuserBR98,BruchhaeuserBR01} at the turn of the last millennium,
cf., e.g., \cite{BruchhaeuserAESW22,BruchhaeuserBR06,BruchhaeuserMB07}.
In order to obtain efficient, adaptive meshes in space and time, the DWR approach
yields an a posteriori error estimator measured in goal quantities of physical 
interest. 
To get such an error estimator, an additional \emph{dual} problem has to be 
solved, which is used for \emph{weighting} the influence of local \emph{residuals} 
on the error.
In this respect, the approximation of these weights is crucial. 
As shown in a comparative numerical study by the authors in \cite{BruchhaeuserBSB19},
in particular considering convection-dominated problems,
the differences between different approximation techniques for the spatial weights
can be significant regarding accuracy and efficiency reasons.
In this work, the focus is now on the numerical investigation of different 
approximation approaches for the temporal weights.
Thus, we consider the following time-dependent convection-diffusion-reaction 
equation given by 
\begin{equation}
\label{eq:1:CDRoriginal}
\begin{array}{r@{\;}c@{\;}l@{\;} @{\,\,}l @{\,\,}l @{\,}l}
\partial_t u
- \nabla \cdot (\varepsilon \nabla u) 
+ \boldsymbol{b} \cdot \nabla u 
+ \alpha u & = & f & \mbox{in } & Q & = \Omega \times I \,, 
\\[1ex]
u(0) & = & u_0 & \mbox{on } & \Sigma_0 & = \Omega \times \{0\} \,, 
\end{array}
\end{equation}
equipped with appropriate boundary conditions. We denote by $Q$ the 
space-time domain, where $\Omega\subset \mathbb{R}^d$, with $d=2$ or $d=3$, is a 
polygonal or polyhedral bounded domain with Lipschitz boundary $\partial\Omega$ 
and $I=(0,T), 0 < T < \infty$, is a finite time interval. 
Well-posedness of \eqref{eq:1:CDRoriginal} and the existence of a sufficiently 
regular solution, such that all of the arguments and terms used below are 
well-defined, are tacitly assumed without mentioning explicitly all technical 
assumptions about the data and coefficients here; cf., e.g., 
\cite{BruchhaeuserEG21} and \cite[Ch. 4]{BruchhaeuserB22} for the details.
Henceforth, for the sake of simplicity, we deal with homogeneous Dirichlet 
boundary values, $u|_{\partial\Omega}=0$. In our numerical examples, 
we also consider more general boundary conditions, whose incorporation is straightforward, 
cf.~\cite[Rem. 4.14]{BruchhaeuserB22}. 

This work is organized as follows. In Sec.~\ref{sec:2:weak}, we present the weak
form as well as the (stabilized) space-time discretization schemes. 
In Sec.~\ref{sec:3:error}, we derive an a posteriori error representation
for the model problem.   
In Sec.~\ref{sec:4:practical}, we introduce two different approximation techniques
for the temporal weights and present the underlying space-time adaptive algorithm.
Finally, in Sec.~\ref{sec:5:numerical}, the results of our numerical comparison 
study are presented. 
\section{Variational Formulation and Discretization in Space and Time}
\label{sec:2:weak}
Introducing the function space  
$V := 
\big\{
v\in L^{2}\big(I; H_0^1(\Omega)\big)\,\big|\,
\partial_{t}v\in L^{2}\big(I; H^{-1}(\Omega)\big)
\big\}\,,$
the variational formulation of problem \eqref{eq:1:CDRoriginal} reads as follows:
\textit{Find $u \in V$ such that}
\begin{equation}
\label{eq:2:WeakCDRsteady}
A(u)(\varphi) := \displaystyle\int_{I}\big\{(\partial_t u,\varphi)
+a(u)(\varphi)
\big\} \mathrm{d} t 
+ (u(0)-u_0\,,\varphi(0))
=
\displaystyle\int_I(f,\varphi)\;\mathrm{d}t
=:F(\varphi) 
\end{equation}
\textit{for all $\varphi \in V$, 
with an inner bilinear form $a(u)(\varphi) := (\varepsilon\nabla u, \nabla \varphi) 
+(\boldsymbol{b}\cdot \nabla u, \varphi) + (\alpha u,\varphi)\,.$
} 
For the discretization in time we use a discontinuous Galerkin method dG($r$) 
with an arbitrary polynomial degree $r\ge 0$.
Let $0=:t_0<t_1<\dots<t_N:=T$ be a set of time points for a partition of 
$I$ into left-open subintervals $I_n:=(t_{n-1},t_n]$, $n=1,\dots,N$,
with time step sizes $\tau_n=t_n-t_{n-1}$ and global discretization
parameter $\tau=\max_{n}\tau_{n}$. Thus, the space-time domain $Q$
is separated into a partition of tensor-product slabs $Q=\Omega\times I_n$.
Let $V_{\tau}^{r}:=
 \big\{u_{\tau}\in L^{2}(I; H_0^1(\Omega))\big|
 u_{\tau}|_{I_{n}}\in \mathcal{P}_{r}(I_{n}; H_0^1(\Omega)),
 u_{\tau}(0)\in L^2(\Omega)\big\}\,,$
where $\mathcal{P}_{r}(I_{n}; H_0^1(\Omega))$ is the space of
polynomials up to degree $r$ on $I_n$ with values in $H_0^1(\Omega)\,.$
Then, the semi-discrete in time scheme reads as follows:
\textit{Find $u_\tau \in V_{\tau}^{r}$ such that
}
\begin{equation}
\label{eq:3:A_tau_u_phi_eq_F_phi}
A_{\tau}(u_\tau)(\varphi_\tau) + (u_{\tau,0}^+ - u_{0},\varphi_{\tau,0}^+)
=
F(\varphi_\tau) 
\quad \forall \varphi_\tau \in V_{\tau}^{r}\,,
\end{equation}
\textit{with $F(\cdot)$ being defined in \eqref{eq:2:WeakCDRsteady} and the 
semi-discrete bilinear form $A_{\tau}(\cdot)(\cdot)$ given by}
\begin{equation}
\label{eq:4:Def_A_tau_u_phi_and_Def_F_phi_tau}
A_{\tau}(u_\tau)(\varphi_\tau)
 := 
\displaystyle\sum_{n=1}^{N}\int_{I_n}\big\{(\partial_{t} u_\tau,\varphi_\tau)
+ a(u_\tau)(\varphi_\tau)
\big\} \mathrm{d} t\; 
+ \displaystyle\sum_{n=2}^{N}([u_\tau]_{n-1},\varphi_{\tau,n-1}^+ )\,.
\end{equation}
Here, the jump $[\cdot]_{n}$ at a time point $t_n$ for some discontinuous in time 
function $u_{\tau}\in V_{\tau}^{r}$ is defined by $[u_{\tau}]_{n} := u_{\tau,n}^{+}
-u_{\tau,n}^{-}\,,$ with $u_{\tau,n}^{\pm}:= \lim_{t\mapsto t_n\pm0} u_\tau(t)$.

\noindent
For the discretization in space, we use a standard continuous Galerkin method 
cG($p$) with an arbitrary polynomial degree $p\geq 1$. For this, we consider a
decomposition $\mathcal{T}_{h}$ of the domain $\Omega$ into disjoint elements
$K$, such that $\overline{\Omega}=\cup_{K\in\mathcal{T}_{h}}\overline{K}$.
Here, we choose the elements $K\in\mathcal{T}_{h}$ to be quadrilaterals
($d=2$) or hexahedrals ($d=3$).
We denote by $h_{K}$ the diameter of $K$ and the global space
discretization parameter is $h:=\max_{K\in\mathcal{T}_{h}}h_{K}$.
Our mesh adaptation process yields locally refined cells, which is
enabled by using hanging nodes; cf.~\cite{BruchhaeuserBR03} for further details. 
Let $V_{\tau h}^{r,p} := \big\{ u_{\tau h}\in V_{\tau}^{r} \big|
u_{\tau h}|_{I_n} \in \mathcal{P}_r(I_n;V_h^{p,n})\,, u_{\tau h}(0) \in V_h^{p,0}
\big\}$ 
be the fully discrete function space.
Now, the fully discrete scheme is given by replacing the semi-discrete functions
in \eqref{eq:3:A_tau_u_phi_eq_F_phi} by the fully discrete ones.
In the implementation we are using tensor-product polynomials on a single slab,
cf. \cite{BruchhaeuserBKB22} for the details.
Since we focus on convection-dominated problems here, the finite element 
approximation needs to be stabilized in order to reduce spurious and non-physical 
oscillations arising close to layers. 
Here, we use the streamline upwind Petrov-Galerkin (SUPG) method \cite{BruchhaeuserHB79,BruchhaeuserBH82}. 
Finally, the stabilized fully discrete scheme reads as follows:
\textit{Find $u_{\tau h} \in V_{\tau h}^{r,p}$ such that}
\begin{equation}
\label{eq:5:stabilized_scheme}
A_{S}(u_{\tau h})(\varphi_{\tau h}) + (u_{\tau h,0}^+ - u_{0},\varphi_{\tau h,0}^+) = 
F(\varphi_{\tau h}) \quad
\forall \varphi_{\tau h} \in V_{\tau h}^{r,p}\,,
\end{equation}
\textit{where} $A_{S}(u_{\tau h})(\varphi_{\tau h}):=
A_{\tau}(u_{\tau h})(\varphi_{\tau h})+S_{A}(u_{\tau h})(\varphi_{\tau h})$, 
\textit{with} $A_{\tau}(\cdot)(\cdot)$ 
\textit{being defined in \eqref{eq:4:Def_A_tau_u_phi_and_Def_F_phi_tau}}
\textit{and} $S_A(\cdot)(\cdot)$ \textit{are the SUPG stabilization terms being 
defined in the Appendix}.
\section{DWR-Based A Posteriori Error Estimation}
\label{sec:3:error}
In this section, we present a DWR-based a posteriori error representation for the 
stabilized problem \eqref{eq:5:stabilized_scheme}.
This representation is given in terms of a user chosen 
goal functional, in general given as
$
J(u)=\int_0^T J_1(u(t))\;\mathrm{d}t + J_2(u(T))\,,
$
where $J_1 \in L^2(I;H^{-1}(\Omega))$ and $J_2 \in H^{-1}(\Omega)$
are three times differentiable functionals and each of them may be zero; 
cf., e.g., \cite{BruchhaeuserR17}.
In order to keep this work  short, we restrict this section to the main 
result only and refer to \cite[Ch. 4.2]{BruchhaeuserB22} for a detailed 
derivation.
\begin{theorem}[Temporal and Spatial Error Representation]
\label{thm:1:ErrorRepresentationDG}
~\\
Let $\{u,z\} \in V \times V$,
$\{u_{\tau},z_{\tau}\} \in V_{\tau}^{r} 
\times V_{\tau}^{r}$,
and
$\{u_{\tau h},z_{\tau h}\} \in V_{\tau h}^{r,p} 
\times V_{\tau h}^{r,p}$
be stationary points of Lagrangian functionals
$\mathcal{L}, \mathcal{L}_{\tau}$, and $\mathcal{L}_{\tau h}$
on different discretization levels
\begin{displaymath}
\begin{array}{r@{\;}c@{\;}l@{\;}l@{\;}}
\mathcal{L}^{\prime}(u,z)(\delta u, \delta z)
= \mathcal{L}_{\tau}^{\prime}(u,z)(\delta u, \delta z)
& = & 0 
& \quad
\forall \, \{\delta u,\delta z\}\in V \times V\,,
\\
\mathcal{L}_{\tau}^{\prime}(u_{\tau},z_{\tau})
(\delta u_{\tau}, \delta z_{\tau})
& = & 0
& \quad 
\forall \, \{\delta u_{\tau},\delta z_{\tau}\}
\in V_{\tau}^{r} \times V_{\tau}^{r}\,,
\\
\mathcal{L}_{\tau h}^{\prime}(u_{\tau h},z_{\tau h})
(\delta u_{\tau h}, \delta z_{\tau h})
& = & 0
& \quad
\forall \, \{\delta u_{\tau h},\delta z_{\tau h}\}
\in V_{\tau h}^{r,p} \times V_{\tau h}^{r,p}\,.
\end{array}
\end{displaymath}
Then, there holds the error representation formulas in space and time, respectively,
\begin{equation}
\label{eq:6:DGreducedJ_u_minus_J_u_tau}
 \begin{array}{r@{\;}c@{\;}l@{\;}}
J(u)-J(u_{\tau}) & = &
\rho_{\tau}(u_{\tau})(z-\tilde{z}_{\tau})
+ \frac{1}{2}\Delta \rho_{\tau}
+ \mathcal{R}_{\tau}\,,
\\[0.5ex]
J(u_{\tau})-J(u_{\tau h}) & = &
\rho_{\tau}(u_{\tau h})(z_{\tau}-\tilde{z}_{\tau h})
+ S_A(u_{\tau h})(\tilde{z}_{\tau h})
+ \frac{1}{2}\Delta \rho_{h} 
+ \mathcal{R}_{h}\,,
\end{array}
\end{equation}
where $\rho_\tau$ denotes the primal residual based on the semi-discrete schemes
\eqref{eq:3:A_tau_u_phi_eq_F_phi}.
Here, $\tilde{z}_{\tau} \in V_{\tau}^{r}$ and $\tilde{z}_{\tau h} \in V_{\tau h}^{r,p}$
can be chosen arbitrarily and $\mathcal{R}_{\tau}, \mathcal{R}_{h}$ are higher-order 
remainder terms with respect to the errors
$u-u_\tau,z-z_\tau$ and
$u_{\tau}-u_{\tau h},z_{\tau}-z_{\tau h}$. 
The explicit presentations of $\mathcal{L}, \mathcal{L}_{\tau}, \mathcal{L}_{\tau h}$ 
and $\rho_\tau$ are given in the Appendix.
\end{theorem}
\begin{mproof}
The proof includes a result regarding the relation between the primal 
and dual residual that goes back to \cite[Prop. 2.3]{BruchhaeuserBR01}.
A detailed derivation of this relation for the convection-diffusion equation can
be found in \cite[Ch. 4.2.2]{BruchhaeuserB22}.
\end{mproof}

\section{Practical Aspects}
\label{sec:4:practical}
The error representation formulas derived in the previous section lead to a 
posteriori error estimators in space and time, whereby the latter serve as 
indicators for  the adaptive mesh refinement process.
Therefore, we present a localized form of the reduced 
error representation formulas. 
Neglecting the higher-order remainder terms in Thm.~\ref{thm:1:ErrorRepresentationDG},
we get the following error representation formulas in space and time.
\begin{equation}
\label{eq:7:DGRemainderlessreducedJ_minus_J}
\begin{array}{r@{\,}c@{\,}l@{\,}c@{\,}l@{\,}c@{\,}l}
J(u)-J(u_{\tau}) & \approx &
\rho_{\tau}(u_{\tau})(z-\tilde{z}_{\tau})
&=:&\eta_{\tau}&=&\displaystyle\sum_{n=1}^{N}\eta_{\tau}^{n}\,,
\\[1.5ex]
J(u_{\tau})-J(u_{\tau h}) & \approx &
\rho_{\tau}(u_{\tau h})(z_{\tau}-\tilde{z}_{\tau h})
+ S_A(u_{\tau h})(\tilde{z}_{\tau h})
&=:&\eta_{h}&=&\displaystyle\sum_{n=1}^{N}\eta_{h}^{n}\,,
\end{array}
\end{equation}
where $\eta_{\tau}^{n}$ and $\eta_{h}^{n}$ denote the local error indicators 
with regard to a single slab $Q_n=\mathcal{T}_{h,n}
\times I_{n}\,,n=1,\dots,N\,;$ cf. \cite[Ch. 4.3.2]{BruchhaeuserB22} for a 
detailed representation.
These indicators include unknown solutions as well as temporal
and spatial weights given by $z-\tilde{z}_{\tau}$ and 
$z_{\tau}-\tilde{z}_{\tau h}$, respectively, that have to be approximated by an 
appropriate technique; cf., e.g., \cite{BruchhaeuserBR03,BruchhaeuserBR12,BruchhaeuserB22}.
Due to the results of a comprehensive numerical comparison study with regard to 
steady-state convection-dominated problems in \cite{BruchhaeuserBSB19,BruchhaeuserB22}, 
we use a higher-order finite elements approach for approximating the spatial 
weights in our numerical examples. 
Now, in this work, we put the focus on the approximation of the temporal weights 
by comparing different approximation techniques similar to the steady-state case.
In the following, we will shortly introduce these different approximation techniques,
where we refer to \cite[Ch. 4.3.2]{BruchhaeuserB22} for more details.

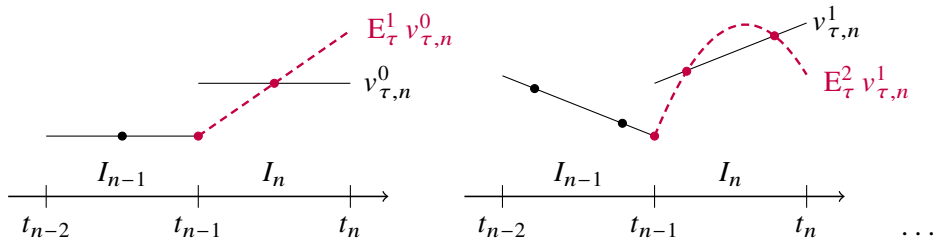
\begin{figure}[hbt!]
\centering
\resizebox{0.88\linewidth}{!}{%
\begin{tikzpicture}
\tikzstyle{ns1}=[line width=0.4]
\tikzstyle{ns2}=[line width=0.15,opacity=.5]
\tikzstyle{ns3}=[line width=0.6]
\tikzstyle{ns4}=[line width=0.9]
\draw[->,ns3] (-8.5,0.0) -- (-3.5,0.0);
\draw[ns1] (-8,-0.15) -- (-8,0.15);
\draw[ns1] (-6,-0.15) -- (-6,0.15);
\draw[ns1] (-4,-0.15) -- (-4,0.15);
\node[] at (-8,-0.4) {$t_{n-2}$};
\node[] at (-7,0.3) {$I_{n-1}$};
\node[] at (-6,-0.4) {$t_{n-1}$};
\node[] at (-5,0.3) {$I_n$};
\node[] at (-4,-0.4) {$t_{n}$};
\draw[domain=-8:-6] plot (\x, {0.8});
\draw[fill=black](-7,0.8)circle(1.5pt);
\draw[HSUred,fill=HSUred](-6,0.8)circle(1.5pt);
\draw[domain=-6:-4] plot (\x, {1.5});
\draw[HSUred,fill=HSUred](-5,1.5)circle(1.5pt);
\draw[HSUred,ns4,densely dashed,domain=-6:-4] plot (\x, {0.7*\x+5});
\node[] at (-3.5,1.5) {$v_{\tau,n}^{0}$};
\node[] at (-3.2,2.2) {\textcolor{HSUred}{$\operatorname{E}_{\tau}^{1}v_{\tau,n}^{0}$}};
\draw[->,ns3] (-2.5,0.0) -- (2.5,0.0);
\draw[ns1] (-2,-0.15) -- (-2,0.15);
\draw[ns1] (0,-0.15) -- (0,0.15);
\draw[ns1] (2,-0.15) -- (2,0.15);
\node[] at (-2,-0.4) {$t_{n-2}$};
\node[] at (-1,0.3) {$I_{n-1}$};
\node[] at (0,-0.4) {$t_{n-1}$};
\node[] at (1,0.3) {$I_n$};
\node[] at (2,-0.4) {$t_{n}$};
\draw[domain=-2:0] plot (\x, {-.4*\x+.8});
\draw[fill=black](-1.5773,1.43)circle(1.5pt);
\draw[fill=black](-0.4227,0.97)circle(1.5pt);
\draw[HSUred,fill=HSUred](0,0.8)circle(1.5pt);
\draw[domain=0:2] plot (\x, {.4*\x+1.5});
\draw[HSUred,fill=HSUred](0.4227,1.66)circle(1.5pt);
\draw[HSUred,fill=HSUred](1.5773,2.13)circle(1.5pt);
\draw[HSUred,ns4,densely dashed,domain=0:2] plot (\x, {-1.0318*\x*\x+2.4708*\x+0.8});
\node[] at (2.4,2.3) {$v_{\tau,n}^{1}$};
\node[] at (2.8,1.5) {\textcolor{HSUred}{$\operatorname{E}_{\tau}^{2}v_{\tau,n}^{1}$}};
\node[] at (3.5,-0.5) {$\dots$};
\end{tikzpicture}
}
\caption{Reconstruction of a discontinuous constant (left) and linear (right) in 
time function using Gauss quadrature points.}
\label{fig:1:hoRetrapolationGauss}
\end{figure}

\begin{flushleft}
\textbf{Approximation by Higher-Order Reconstruction (for short, hoRe approach)} 
\end{flushleft}
\vspace{-0.2cm}
\noindent
For the higher-order reconstruction approach, we introduce the following operator
\begin{displaymath}
z-\tilde{z}_{\tau}  \approx 
\operatorname{E}_{\tau}^{r+1}z_{\tau}-z_{\tau}\,,
\end{displaymath}
where $\operatorname{E}_{\tau}^{r+1}$ denotes a reconstruction in time operator 
that acts on $I_n$ and lifts the solution to a piecewise polynomial of degree 
$(r+1)$ in time, cf. Fig.~\ref{fig:1:hoRetrapolationGauss}. 
%
\vspace{-0.3cm}
%

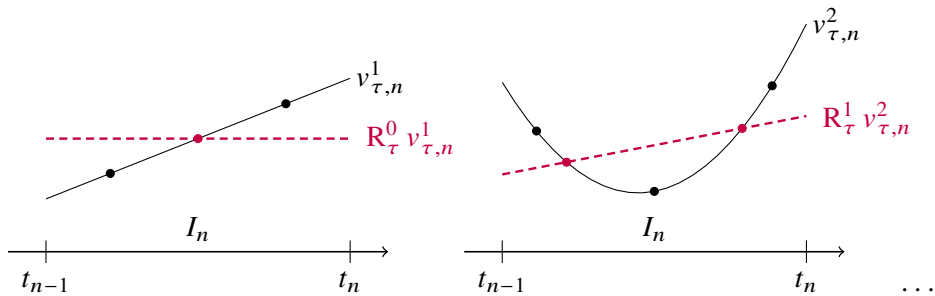
\begin{figure}[hbt!]
\centering
\resizebox{0.88\linewidth}{!}{%
\begin{tikzpicture}
\tikzstyle{ns1}=[line width=0.4]
\tikzstyle{ns2}=[line width=0.15,opacity=.5]
\tikzstyle{ns3}=[line width=0.6]
\tikzstyle{ns4}=[line width=0.9]

\draw[->,ns3] (-2.5,0.0) -- (2.5,0.0);
\draw[ns1] (-2,-0.15) -- (-2,0.15);
\draw[ns1] (2,-0.15) -- (2,0.15);
\node[] at (-2,-0.4) {$t_{n-1}$};
\node[] at (0,0.3) {$I_n$};
\node[] at (2,-0.4) {$t_{n}$};
\draw[domain=-2:2] plot (\x, {.4*\x+1.5});
\draw[fill=black](-1.1547,1.0381)circle(1.5pt);
\draw[HSUred,fill=HSUred](0,1.5)circle(1.5pt);
\draw[fill=black](1.1547,1.9619)circle(1.5pt);
\draw[HSUred,ns4,densely dashed,domain=-2:2] plot (\x, {+1.5});
\node[] at (2.4,2.3) {$v_{\tau,n}^{1}$};
\node[] at (2.8,1.5) {\textcolor{HSUred}{$\operatorname{R}_{\tau}^{0}v_{\tau,n}^{1}$}};
\draw[->,ns3] (3.5,0.0) -- (8.5,0.0);
\draw[ns1] (4,-0.15) -- (4,0.15);
\draw[ns1] (8,-0.15) -- (8,0.15);
\node[] at (4,-0.4) {$t_{n-1}$};
\node[] at (6,0.3) {$I_n$};
\node[] at (8,-0.4) {$t_{n}$};
\draw[domain=4:8] plot (\x, {0.45833*\x*\x-5.3063*\x+16.138});
\draw[fill=black](4.4508,1.6)circle(1.5pt);
\draw[HSUred,fill=HSUred](4.8453,1.1875)circle(1.5pt);
\draw[fill=black](6,0.8)circle(1.5pt);
\draw[HSUred,fill=HSUred](7.1547,1.6348)circle(1.5pt);
\draw[fill=black](7.5492,2.2)circle(1.5pt);
\draw[HSUred,ns4,densely dashed,domain=4:8] plot (\x, {0.1937*\x+0.25});
\node[] at (8.4,3.0) {$v_{\tau,n}^{2}$};
\node[] at (8.8,1.76) {\textcolor{HSUred}{$\operatorname{R}_{\tau}^{1}v_{\tau,n}^{2}$}}; 
\node[] at (9.5,-0.5) {$\dots$};
\end{tikzpicture}
}
\caption{Restriction of a discontinuous linear (left) and quadratic (right) in 
time function using Gauss quadrature points.}
\label{fig:2:hoRestriction}
\end{figure}

\vspace{-0.4cm}
\begin{flushleft}
\textbf{Approximation by Higher-Order Finite Elements (for short, hoFE approach)}
\end{flushleft}
\vspace{-0.2cm}
\noindent
This approach aims to increase the influence of the weights by approximating the
dual solution using higher-order finite elements, i.e. $z_{\tau}$ is computed in 
$V_{\tau}^{s}, s\geq r+1$. 
Therefore, we introduce the following restriction operator in time given by 
\begin{displaymath}
z-\tilde{z}_{\tau} \approx
z_{\tau}-\operatorname{R}_{\tau}^{r}z_{\tau}\,,
\end{displaymath}
where $\operatorname{R}_{\tau}^{r}$ acts on $I_n$ and restricts the solution to
a polynomial of degree $r < s$; cf. Fig.~\ref{fig:2:hoRestriction}.
Note that this approximation is valid, since the quantity $\tilde{z}_{\tau}$ can 
be chosen arbitrarily in the finite element space 
$V_{\tau}^{s}$.
%
\vspace{-0.2cm}
%

\noindent
Now, the last step in making the error indicators computable is to replace all 
unknown solutions in Eqs.~\eqref{eq:7:DGRemainderlessreducedJ_minus_J} by the 
computed fully discrete solutions, i.e. we approximate
\begin{equation}
\label{eq:8:ApproxEE}
\begin{array}{r@{\,}c@{\,}l@{\,}}
\eta_{\tau}&\approx&
\rho_{\tau}(u_{\tau h})(\operatorname{E}_{\tau}^{r+1}z_{\tau h}-z_{\tau h})
=:\tilde{\eta}_{\tau}^{\textnormal{hoRe}}\,, 
\textnormal{ or}
\\[0.5ex]
\eta_{\tau}&\approx&
\rho_{\tau}(u_{\tau h})(z_{\tau h}-\operatorname{R}_{\tau}^{r}z_{\tau h})
=:\tilde{\eta}_{\tau}^{\textnormal{hoFE}}\,,
\textnormal{ and}
\\[0.5ex]
\eta_{h}&\approx&
\rho_{\tau}(u_{\tau h})(z_{\tau h}-\operatorname{R}_{h}^{p}z_{\tau h})
+ S_A(u_{\tau h})(\operatorname{R}_{h}^{p}z_{\tau h})
=:\tilde{\eta}_{h}\,,
\end{array}
\end{equation}
where $\operatorname{R}_{h}^{p}$ denotes a restriction in space operator that 
acts on a spatial cell and restricts the solution to a polynomial of degree
$p<q$ on the corresponding reference element; cf., e.g., \cite{BruchhaeuserBR03,BruchhaeuserB22}.
Finally, we present a survey-like version of our adaptive solution algorithm 
using tensor-product space-time slabs $Q_n$, where we refer to 
\cite[Ch. 4.3.3]{BruchhaeuserB22} and \cite{BruchhaeuserBKB22} for a more 
detailed version and further information.

\vskip2ex
\begin{flushleft}
\textbf{Algorithm: Goal-Oriented Space-Time Adaptivity}
\end{flushleft}

1. Set $\ell=1$ and generate initial space-time slabs
$Q_n^1=\mathcal{T}_{h,n}^1\times I_n^1\,, n=1,\dots,N^1$.
\vskip1ex
2. Compute $\boldsymbol{u_{\tau h}^\ell}$ by solving the \textbf{stabilized primal} 
problem \eqref{eq:5:stabilized_scheme}.
\vskip1ex
3. Compute $\boldsymbol{z_{\tau h}^\ell}$ by solving the \textbf{stabilized dual} 
problem \eqref{eq:12:DualProblems}.
\vskip1ex
4. Evaluate the a posteriori space-time \textbf{error indicators} 
$\boldsymbol{\eta_\tau^\ell}$ and $\boldsymbol{\eta_h^\ell}$.
\vskip1ex
\quad\textbf{if} $|\eta_\tau^\ell|>\omega|\eta_h^\ell|\,,1.5\leq\omega\leq3.5$: 
Adapt the \textbf{temporal} mesh only.
\vskip1ex
\quad\textbf{else if} $|\eta_h^\ell|>\omega|\eta_\tau^\ell|$: 
Adapt the \textbf{spatial} mesh only.
\vskip1ex
\quad\textbf{else}: Adapt both the \textbf{temporal} as well as the 
\textbf{spatial} mesh.
\vskip1ex
5. Increase $\ell$ to $\ell +1$ and return to step~2.


\section{Numerical Examples}
\label{sec:5:numerical}
In this section, we investigate our space-time adaptive algorithm with regard to
the different approximation techniques of the temporal weights introduced in the
previous section. 
Thereby, we study accuracy, efficiency as well as stability properties 
by means of two numerical examples.
The first one is an academic test problem including some challenging behavior
in time, whereas the second one deals with a benchmark problem of 
convection-dominated transport problems.
In order to focus on the results, we restrict the description of the test settings
to a minimum and refer to \cite{BruchhaeuserB22} for the details.
For measuring the accuracy of the error estimator, we consider an effectivity 
index $\mathcal{I}_{\textnormal{eff}}:= \left|\frac{
\tilde{\eta}_{\tau}
+
\tilde{\eta}_{h}}
{J(u)-J(u_{\tau h})}\right|$.
Moreover, in the second example the SUPG stabilization parameter 
is given by
$\delta_K = \delta_0 \cdot h_K \,, \quad 0.1\leq\delta_0\leq1\,,$
where $h_K$ denotes the cell diameter of the spatial mesh cell $K$, 
cf.~Rem.~\ref{rem:1:supg}.

\paragraph{Example 1: Rotating Hill with Changing Orientation} %

As first example, we consider a counterclockwise rotating hill on the unit square
$\Omega=(0,1)^2$, which additionally changes its height and orientation over the 
period $T=1$. 
The corresponding exact solution is given by (cf. \cite[Ex. 4.2]{BruchhaeuserB22}) 
\begin{equation}
\label{eq:9:uexactKB2}
u(\boldsymbol x, t) := \frac{1 + a \cdot ( (x_1 - m_1(t))^2 + (x_2 - m_2(t))^2  ) }{\nu_1(t) \cdot s \cdot \arctan( \nu_2(t) )}
\end{equation}
with $m_1(t) := \frac{1}{2} + \frac{1}{4} \cos(2 \pi t)$ and
$m_2(t) := \frac{1}{2} + \frac{1}{4} \sin(2 \pi t)$, and,
$\nu_1(\hat{t}) := -1, \nu_2(\hat{t}) := 5 \pi \cdot(4\hat{t}-1),$ 
for $\hat{t} \in [0,0.5)$ and 
$\nu_1(\hat{t}) := 1, \nu_2(\hat{t}) := 5\pi\cdot(4(\hat{t}-0.5)-1),$ 
for $\hat{t} \in [0.5,1)$, $\hat{t} = t-k, k \in \mathbf{N}_0$, and,
scalars $s=\frac{1}{3}, a_0 = 50$.
We choose $\varepsilon = 1$, $\boldsymbol{b} = (2,3)^{\top}$ and $\alpha = 1$. 
Furthermore, no stabilization ($\delta_K = 0$) is used in this test case.
The goal quantity is chosen to control the global $L^2(L^2)$-error in space 
and time of $e:=u-u_{\tau h}$, given by
$J_Q(u)= \frac{1}{\|e\|_{(0,T)\times\Omega}}\int_I(u,e)\mathrm{d}t$.
In Tables~\ref{tab:1:KB2-hoReGauss}--\ref{tab:2:KB2-hoFE}, we present the 
results obtained by the two different approximation techniques. 
Here, $N$ denotes the total number of slabs, $N_{K}^{\text{max}}$ the number of 
cells on the finest spatial mesh within the current loop, 
and $N_{\text{DoF}}^{\text{tot}}$ the total space-time degrees of freedom.
Further, $e^{p,r,q,s}$ or rather a cG($p$)-dG($r$)/cG($q$)-dG($s$)
discretization corresponds to a primal solution approximation $u_{\tau h}^{p,r}$ 
in cG($p$)-dG($r$) and a dual solution approximation $z_{\tau h}^{q,s}$ in
cG($q$)-dG($s$).

Regarding the accuracy of the underlying error estimator, as given by the last 
column of Tables~\ref{tab:1:KB2-hoReGauss}--\ref{tab:2:KB2-hoFE}, it can be 
summarized that both approaches show a good quantitative estimation of the 
discretization error as the respective effectivity index $\mathcal{I}_{\mathrm{eff}}$
is around one for all underlying discretization cases. 
We note slightly better results using the hoRE approach regarding accuracy and 
stability of the underlying estimators. 
With regard to efficiency reasons it is essential to ensure an equilibrated 
reduction of the temporal as well as spatial discretization error. 
Referring to this, we point out a good equilibration of the spatial and temporal
error indicators $\tilde{\eta}_{h}$ and $\tilde{\eta}_{\tau}$ in the course of the 
refinement process (columns six and seven of 
Tables~\ref{tab:1:KB2-hoReGauss}--\ref{tab:2:KB2-hoFE}).

\begin{table}[!h]
\centering
\resizebox{\linewidth}{!}{%
\begin{tabular}{c || rrr | c | cc | c || rrr | c | cc | c} 
\hline
 & \multicolumn{7}{c}{cG($1$)-dG($0$)/cG($2$)-dG($0$)} ||& \multicolumn{7}{c}{cG($2$)-dG($1$)/cG($3$)-dG($1$)}
\\ 
\hline
$\ell$ & $N$ & $N_{K}^{\text{max}}$ & $N_{\text{DoF}}^{\text{tot}}$ & 
$J(e^{1,0,2,0})$ & $\tilde{\eta}_h$ & 
$\tilde{\eta}_{\tau}^{\textnormal{hoRe}}$ & 
$\mathcal{I}_{\mathrm{eff}}$ 
& $N$ & $N_{K}^{\text{max}}$ & $N_{\text{DoF}}^{\text{tot}}$ & 
$J(e^{1,1,2,1})$
& $\tilde{\eta}_h$ & 
$\tilde{\eta}_{\tau}^{\textnormal{hoRe}}$ & 
$\mathcal{I}_{\mathrm{eff}}$
\\ 
\hline
1 & 25 &  16 &   625 & 2.74e-02 & 1.89e-02 & 5.11e-03 & 0.88 & 20 &  16 &  1000 & 2.49e-02 & 2.55e-02 & 5.54e-04 & 1.04\\ 
2 & 25 &  40 &  1067 & 1.57e-02 & 6.72e-03 & 4.58e-03 & 0.72 & 20 &  40 &  1724 & 1.19e-02 & 1.08e-02 & 1.15e-03 & 1.01\\ 
3 & 31 & 100 &  3039 & 8.74e-03 & 3.08e-03 & 7.36e-03 & 1.19 & 20 &  88 &  3404 & 7.83e-03 & 6.81e-03 & 1.55e-03 & 1.06\\ 
4 & 37 & 100 &  3583 & 7.18e-03 & 3.73e-03 & 8.19e-03 & 1.65 & 20 & 148 &  5352 & 5.64e-03 & 3.34e-03 & 2.01e-03 & 0.95\\ 
5 & 44 & 100 &  4210 & 6.92e-03 & 3.87e-03 & 5.62e-03 & 1.37 & 20 & 256 &  7600 & 4.83e-03 & 1.66e-03 & 2.54e-03 & 0.87\\ 
6 & 52 & 160 &  7408 & 4.66e-03 & 1.85e-03 & 3.73e-03 & 1.19 & 37 & 256 & 14266 & 2.57e-03 & 2.40e-03 & 3.98e-04 & 1.08\\ 
7 & 62 & 160 &  8864 & 4.56e-03 & 1.87e-03 & 2.76e-03 & 1.01 & 37 & 424 & 19190 & 2.09e-03 & 1.70e-03 & 4.96e-04 & 1.04\\   
8 & 88 & 208 & 15276 & 3.24e-03 & 1.19e-03 & 2.20e-03 & 1.04 & 37 & 724 & 29182 & 1.84e-03 & 1.01e-03 & 5.53e-04 & 0.85\\ 
\hline
\end{tabular}
}
\caption{Adaptive refinement including effectivity indices with goal quantity 
$J_Q$, $\varepsilon=1$, $\delta_0=0$, and $\omega=1.5$ for Example~1,
using a \textbf{higher-order reconstruction} strategy for the temporal weights.
}
\label{tab:1:KB2-hoReGauss} 
\end{table}
\begin{table}[!h]
\centering
\resizebox{\linewidth}{!}{%
\begin{tabular}{c || rrr | c | cc | c || rrr | c | cc | c} 
\hline
 & \multicolumn{7}{c}{cG($1$)-dG($0$)/cG($2$)-dG($1$)} ||& \multicolumn{7}{c}{cG($1$)-dG($1$)/cG($2$)-dG($2$)}
\\ 
\hline
$\ell$ & $N$ & $N_{K}^{\text{max}}$ & $N_{\text{DoF}}^{\text{tot}}$ & 
$J(e^{1,0,2,1})$ & $\tilde{\eta}_h$ & 
$\tilde{\eta}_{\tau}^{\textnormal{hoRe}}$ & 
$\mathcal{I}_{\mathrm{eff}}$ 
& $N$ & $N_{K}^{\text{max}}$ & $N_{\text{DoF}}^{\text{tot}}$ & 
$j(e^{1,1,2,2})$
& $\tilde{\eta}_h$ & 
$\tilde{\eta}_{\tau}^{\textnormal{hoRe}}$ & 
$\mathcal{I}_{\mathrm{eff}}$
\\ 
\hline
1 &  25 &  16 &   625 & 2.74e-02 & 2.18e-02 & 1.10e-02 & 1.19 & 20 &  16 &  1000 & 2.49e-02 & 2.63e-02 & 1.27e-03 & 1.10\\ 
2 &  25 &  40 &  1053 & 1.58e-02 & 8.98e-03 & 1.89e-02 & 1.77 & 20 &  40 &  1808 & 1.12e-02 & 1.05e-02 & 2.86e-03 & 1.18\\ 
3 &  33 &  40 &  1409 & 1.30e-02 & 9.92e-03 & 1.04e-02 & 1.56 & 24 & 160 &  7796 & 2.80e-03 & 2.79e-03 & 2.05e-03 & 1.72\\ 
4 &  42 &  88 &  3842 & 7.38e-03 & 5.15e-03 & 1.02e-02 & 2.07 & 28 & 316 & 14340 & 1.58e-03 & 7.02e-04 & 1.44e-03 & 1.35\\ 
5 &  54 &  88 &  4854 & 6.64e-03 & 5.32e-03 & 6.55e-03 & 1.78 & 33 & 316 & 16990 & 1.52e-03 & 7.59e-04 & 1.07e-03 & 1.19\\ 
6 &  70 & 160 &  8788 & 4.16e-03 & 3.12e-03 & 5.54e-03 & 2.08 & 46 & 544 & 39612 & 8.79e-04 & 4.08e-04 & 5.04e-04 & 1.03\\ 
7 &  91 & 160 & 11395 & 4.03e-03 & 3.06e-03 & 4.07e-03 & 1.76 & 55 & 820 & 74310 & 5.15e-04 & 2.02e-04 & 3.77e-04 & 1.12\\ 
8 & 105 & 160 & 16357 & 3.16e-03 & 2.15e-03 & 3.78e-03 & 1.87 & 66 & 820 & 89752 & 4.87e-04 & 2.03e-04 & 2.56e-04 & 0.94\\ 
\hline
\end{tabular}
}
\caption{Adaptive refinement including effectivity indices with goal quantity 
$J_Q$, $\varepsilon=1$, $\delta_0=0$, and $\omega=1.5$ 
for Example~1, using a \textbf{higher-order finite elements} strategy for the 
temporal weights.}
\label{tab:2:KB2-hoFE} 
\end{table}

\paragraph{Example 2: Moving Hump with Circular Layer} %

The second example is a benchmark for convection-dominated problems that goes back 
to \cite{BruchhaeuserJS08}.
The (analytical) solution is characterized by a hump changing its height in the 
course of the time, given by (cf. \cite[Ex. 4.3]{BruchhaeuserB22}) 
\begin{equation}
\label{eq:10:uexactMH}
\begin{array}{r@{\;}c@{\;}l@{\;}}
u(\boldsymbol{x}, t) & := &
\frac{16}{\pi}\sin(\pi t)x_1(1-x_1)x_2(1-x_2)
\\
& & \cdot
\big( \frac{1}{2} + 
\arctan\big[2\varepsilon^{-\frac{1}{2}}\big(
r_0^2-(x_1-x_1^0)^2-(x_2-x_2^0)^2\big)
\big]
\big)
\end{array}
\end{equation}
The problem is defined on $\Omega\times I := (0,1)^2\times (0,0.5]$ with the 
scalars $r_0=0.25$ and $x_1^0=x_2^0=0.5$.
We choose $\boldsymbol{b} = (2,3)^{\top}$ and $\alpha = 1$. 
The goal quantity is chosen to control the $L^2$-error
$e_N^-=u(\boldsymbol{x},T)-u_{\tau h}(T^{-})$, at the final time point $T$, 
given by $J_T(u)= \big(u(\boldsymbol{x},T),e_N^-\big)/\|e_N^-\|$.
%
\begin{figure}[hbt!]
\centering
\begin{minipage}{.32\linewidth}
\centering
\begin{tikzpicture}
\begin{axis}[%
width=1.5in,
height=2.0in,
scale only axis,
/pgf/number format/.cd, 1000 sep={},
xmode=log,
ymin=0.,
ymax=3,
yminorticks=true,
legend style={legend pos=north west},
legend style={draw=black,fill=white,legend cell align=left,font=\normalsize},
legend entries = {
  {$\{u_{\tau h}^{1,0}, z_{\tau h}^{2,0}\}$ hoRe},
  {$\{u_{\tau h}^{1,0}, z_{\tau h}^{2,1}\}$ hoFE},
}
]
\node at (axis cs:3200, 1.95) {\fbox{\small$\varepsilon=10^{-3}$}};
\addplot [
color=red,
densely dashdotted,
line width=1.0pt,
mark=square*,
mark size = 1.,
mark options={solid,red}
]
table[row sep=crcr]{
810  0.7917\\  
1470  1.1436\\  
2274  1.4567\\ 
5057  1.4204\\ 
8433  1.6196\\ 
16491  1.5311\\ 
17869  1.4151\\ 
23351  1.8267\\ 
25978  1.6725\\ 
38944  1.6436\\ 
49438  1.7574\\ 
71522  1.8837\\ 
78600  1.5748\\ 
113430  1.4958\\ 
157064  1.4240\\  
254340  1.4063\\ 
371286  1.4649\\ 
425193  1.2665\\ 
558126  1.1518\\  
850679  1.0031\\ 
1121719  0.9914\\ 
1229939  0.9450\\ 
2250302  0.9001\\ 
};
\addplot [
color=blue,
densely dotted,
line width=1.0pt,
mark=x,
mark options={solid,blue}
]
table[row sep=crcr]{
810  1.0930\\ 
1644  1.1015\\ 
2938  1.0488\\ 
5204  1.0205\\ 
13538  1.0487\\ 
28542  1.1028\\ 
67692  1.0983\\ 
104276  1.0814\\ 
257599  1.0595\\ 
399246  0.9788\\ 
1013318  0.9895\\ 
1594792  1.0427\\ 
2269998  0.9266\\ 
};
\end{axis}
\end{tikzpicture}
\end{minipage}
\begin{minipage}{.32\linewidth}
\centering
\begin{tikzpicture}
\begin{axis}[%
width=1.5in,
height=2.0in,
scale only axis,
/pgf/number format/.cd, 1000 sep={},
xmode=log,
ymin=0.,
ymax=3,
yminorticks=true,
legend style={legend pos=north west},
legend style={draw=black,fill=white,legend cell align=left,font=\normalsize},
legend entries = {
  {$\{u_{\tau h}^{1,0}, z_{\tau h}^{2,0}\}$ hoRe},
  {$\{u_{\tau h}^{1,0}, z_{\tau h}^{2,1}\}$ hoFE},
}
]
\node at (axis cs:3100, 1.95) {\fbox{\small$\varepsilon=10^{-4}$}};
\addplot [
color=red,
densely dashdotted,
line width=1.0pt,
mark=square*,
mark size = 1.,
mark options={solid,red}
]
table[row sep=crcr]{
810     0.4880\\ 
1570    0.5969\\  
2670    0.8421\\  
5890    1.0504\\  
15752   0.9925\\ 
34330   1.2999\\ 
65247   1.2824\\ 
141846  1.5532\\ 
224842  1.0572\\  
479340  1.1552\\ 
981776  0.9325\\ 
1745235 0.8445\\ 
};
\addplot [
color=blue,
densely dotted,
line width=1.0pt,
mark=x,
mark options={solid,blue}
]
table[row sep=crcr]{
810  0.7067\\  
1388  0.8749\\ 
2458  0.8292\\  
3896  0.9173\\ 
7124  0.8055\\ 
10471  0.9884\\ 
13570  1.0019\\  
22500  1.1391\\ 
41684  1.1464\\ 
60283  1.1367\\  
105613  1.2249\\ 
168377  1.1262\\ 
239528  1.0963\\ 
362271  1.0302\\ 
576010  0.9444\\ 
1148692  0.8933\\ 
1572170  0.9269\\ 
};
\end{axis}
\end{tikzpicture}
\end{minipage}
\begin{minipage}{.32\linewidth}
\centering
\begin{tikzpicture}
\begin{axis}[%
width=1.5in,
height=2.0in,
scale only axis,
/pgf/number format/.cd, 1000 sep={},
xmode=log,
ymin=0.,
ymax=3.,
yminorticks=true,
legend style={legend pos=north west},
legend style={draw=black,fill=white,legend cell align=left,font=\normalsize},
legend entries = {
  {$\{u_{\tau h}^{1,0}, z_{\tau h}^{2,0}\}$ hoRe},
  {$\{u_{\tau h}^{1,0}, z_{\tau h}^{2,1}\}$ hoFE},
}
]
\node at (axis cs:11500, 1.95) {\fbox{\small$\varepsilon=10^{-6}$}};
\addplot [
color=red,
densely dashdotted,
line width=1.0pt,
mark=square*,
mark size = 1.,
mark options={solid,red}
]
table[row sep=crcr]{
2890  0.0286\\  
4046  0.0104\\  
7051  0.0509\\  
9924  0.0353\\  
17035  0.0918\\  
27067  0.0468\\ 
52772  0.2300\\  
78232  0.2858\\  
119460  0.4196\\  
179478  0.5526\\  
285564  0.4776\\ 
443048  0.6842\\ 
677052  1.0373\\ 
1779166  1.1360\\ 
3492586  1.0632\\ 
4999150 1.0130\\
7672162 1.0316\\
};
\addplot [
color=blue,
densely dotted,
line width=1.0pt,
mark=x,
mark options={solid,blue}
]
table[row sep=crcr]{
2890 0.0849\\ 
3456 0.1219\\ 
6792 0.1353\\ 
13448 0.1505\\ 
23068 0.3638\\ 
43208 0.4877\\ 
61312 0.7229\\ 
117303 0.9193\\  
222226 0.8663\\ 
474275 0.9278\\ 
840502 0.9755\\ 
1731154 1.0139\\ 
3545763 0.9740\\ 
7419449 0.9533\\ 
};
\end{axis}
\end{tikzpicture}
\end{minipage}
\caption{Comparison of effectivity indices $\mathcal{I}_{\mathrm{eff}}$ 
(over total space-time DoFs) regarding both approximation 
approaches for varying diffusion coefficients $\varepsilon$, using 
$\delta_0=1$, $\omega = 2$ and $J_T$ for Example~2.}
\label{fig:4:IeffMH}
\end{figure}
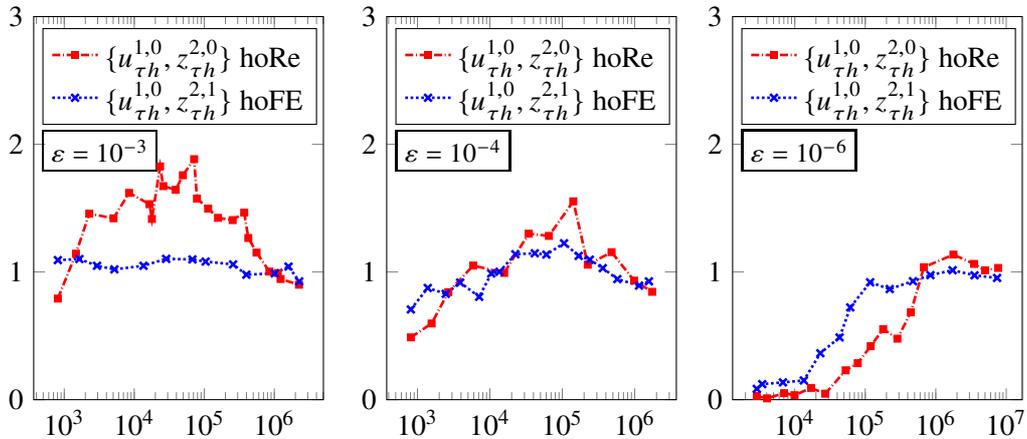
%
In Fig.~\ref{fig:4:IeffMH}, we compare both approximation approaches regarding
the development of the respective effectivity indices for a sequence of 
decreasing diffusion coefficients.
Overall, the the desired value of one is reached in the course of the refinement 
process for both approaches in all considered cases of $\varepsilon$, 
but there are slightly differences with regard to accuracy, efficiency and stability
reasons.  
For $\varepsilon = 10^{-3}, 10^{-4}$, we observe a slight overestimation of the
exact error by using the hoRe approach.
In contrast, the hoFE approach seems to be more stable. 
The more challenging case of $\varepsilon = 10^{-6}$ shows a slightly different 
behavior. 
Here, both approaches show up with small values for the effectivity index at the 
beginning. It takes some refinement processes until this value increases.
Again, the hoFE approach shows slightly better results with regard to accuracy and
efficiency evident by an faster achievement of the desired value of one.
This behavior is in good agreement to the results obtained for the 
steady-state case, even though the differences between the both approaches is not
that significant here, cf.~Example~1 in \cite{BruchhaeuserBSB19} and 
\cite[Ch. 3]{BruchhaeuserB22}. 

\section{Conclusion}
\label{sec:6:conclusion}
In this work, we compared two different approaches for approximating the temporal
weights occurring within the error estimators based on the DWR method. We studied
accuracy, efficiency and stability reasons by means of a convection-dominated
model problem. 
Robustness of the underlying space-time adaptive algorithm was shown and 
effectivity indices close to one were obtained for both approximation techniques.
Even though the higher-order finite elements approach shows slightly better 
results with regard to accuracy and stability (cf.~Fig.~\ref{fig:4:IeffMH}), 
the differences to the higher-order reconstruction approach in all was not that 
significant as for the steady-state case (cf.~\cite{BruchhaeuserBSB19,BruchhaeuserB22}), 
even less if we take into account the higher computational costs for solving the 
dual problem.

\section*{Acknowledgement}
We acknowledge U. K\"ocher for his support in the design and 
implementation of the underlying software \texttt{dwr-condiffrea}; 
cf. the software project \texttt{DTM++.Project/dwr} \cite{BruchhaeuserKBB19}
that is based on the open source finite element library \texttt{deal.II} 
\cite{BruchhaeuserABBFGHHKKMMPPSTWZ21}.

\section*{Appendix}
The SUPG stabilization terms for the stabilized primal problem \eqref{eq:5:stabilized_scheme}
are defined as
\begin{displaymath}
\begin{array}{l@{}l}
S_{A}(u_{\tau h})(\varphi_{\tau h}) &{} :=  \displaystyle\sum_{n=1}^N\int_{I_n} 
\sum\limits_{K\in \mathcal{T}_h}\delta_K\big( r(u_{\tau h}), 
\boldsymbol{b}\cdot\nabla\varphi_{\tau h}\big)_K \,\mathrm{d} t \\[0ex]
+ \displaystyle\sum\limits_{n = 2}^{N}\sum\limits_{K\in\mathcal{T}_h} 
&{}\hspace{-0.2cm} \delta_K
\big(\left[u_{\tau h}\right]_{n-1}, 
\boldsymbol{b}\cdot\nabla\varphi_{\tau h,n-1}^+\big)_{K} 
+ \displaystyle\sum\limits_{K\in\mathcal{T}_h}
\delta_K \big( u_{\tau h,0}^{+} - u_0,
\boldsymbol{b}\cdot\nabla\varphi_{\tau h,0}^{+} \big)_{K}\,,\vspace{0.2cm} \\[0ex]
\hspace{0.8cm}r(u_{\tau h}) &{} :=  \partial_{t} u_{\tau h} 
+ \boldsymbol{b}\cdot\nabla u_{\tau h} 
- \nabla\left(\varepsilon\nabla u_{\tau h}\right)  + \alpha u_{\tau h} - f\,.
\end{array}
\end{displaymath}
\begin{remark}
\label{rem:1:supg}
For time-dependent convection-diffusion problems an optimal error estimate for 
$\delta_K=\mathrm{O}(h_K)$ is derived in \cite{BruchhaeuserJN11}.
\end{remark}
\noindent
The Lagrangian functionals
$\mathcal{L}: V\times V \rightarrow \mathbb{R}\,,$
$\mathcal{L}_\tau: V_{\tau}^{r} 
\times V_{\tau}^{r}\rightarrow \mathbb{R}\,,$
and
$\mathcal{L}_{\tau h}:
V_{\tau h}^{r,p} 
\times V_{\tau h}^{r,p} \rightarrow \mathbb{R}$ 
are defined by
\begin{equation}
\label{eq:11:Lagrangians}
\begin{array}{r@{\;}c@{\;}l@{\;}}
\mathcal{L}(u,z) & := & J(u)
+ F(z)
- A(u)(z)\,,
\\[1ex]
\mathcal{L}_{\tau}(u_\tau,z_\tau) & := &
J(u_{\tau}) + F(z_{\tau})
- A_{\tau}(u_{\tau})(z_{\tau})\,,
\\[1ex]
\mathcal{L}_{\tau h}(u_{\tau h},z_{\tau h}) & := &
J(u_{\tau h})
+ F (z_{\tau h})
- A_{S}(u_{\tau h})(z_{\tau h})
\end{array}
\end{equation}
Here, the Lagrange multipliers $z$, $z_\tau,$ and $z_{\tau h}$
are called dual variables in contrast to the primal variables
$u$, $u_\tau,$ and
$u_{\tau h}$; cf. \cite{BruchhaeuserBR12,BruchhaeuserBR01}.
Considering the directional derivatives of the Lagrangian functionals, also
known as
G\^{a}teaux derivatives, with respect to their first argument, i.e.\
$
\mathcal{L}^{\prime}_{u}(u,z)(\varphi) :=
\lim_{t\neq0,t\rightarrow 0}
t^{-1}\big\{\mathcal{L}(u + t\varphi,z)
-\mathcal{L}(u,z)\big\},
\quad \varphi \in V\,,
$
leads to the  dual problems: 
Find the continuous dual solution $z \in V$, the semi-discrete
dual solution $z_{\tau} \in V_{\tau}^{r}$ and the fully
discrete dual solution $z_{\tau h} \in V_{\tau h}^{r,p}$,
respectively, such that
\begin{equation}
\label{eq:12:DualProblems}
\begin{array}{r@{\,}c@{\,}l@{\,}l@{\,}}
A^{\prime}(u)(\varphi,z)
&=&
J^{\prime}(u)(\varphi)
&
\forall \varphi \in V\,,
\\[1.ex]
A_{\tau}^{\prime}(u_{\tau})(\varphi_{\tau},z_{\tau})
& = &
J^{\prime}(u_{\tau})(\varphi_{\tau})
&
\forall \varphi_{\tau}\in V_{\tau}^{r}\,,
\\[1.ex]
A_{S}^{\prime}(u_{\tau h})(\varphi_{\tau h},z_{\tau h})
& = &
J^{\prime}(u_{\tau h})(\varphi_{\tau h})
& \forall \varphi_{\tau h}\in V_{\tau h}^{r,p}\,,
\end{array}
\end{equation}
where we refer to our work \cite[Ch. 4.2.1]{BruchhaeuserB22} for a detailed description of the
adjoint bilinear forms $A^{\prime},A_{\tau}^{\prime},A_S^{\prime}$ as well as 
the dual right hand side term $J^\prime$. 
Finally, the primal residual based on the semi-discrete scheme is defined by means of the
G\^{a}teaux derivatives of the Lagrangian functionals in the following way:
\begin{displaymath} 
\rho_{\tau}(u)(\varphi) := 
\mathcal{L}_{\tau,z}^{\prime}(u,z)(\varphi)
 =  G_{\tau}(\varphi)-A_{\tau}(u)(\varphi)\,,
\end{displaymath}

\end{document}